\documentclass[12pt]{amsart}
\usepackage{amssymb}
\usepackage[all]{xy}
\usepackage{lscape}

\theoremstyle{plain} 
\newtheorem{theor}[equation]{Theorem}
\newtheorem{cor}[equation]{Corollary}
\newtheorem{lem}[equation]{Lemma}
\newtheorem{conjecture}[equation]{Conjecture}
\newtheorem{proposition}[equation]{Proposition}

\theoremstyle{definition}
\newtheorem{defin}[equation]{Definition}

\theoremstyle{remark}
\newtheorem{rem}[equation]{Remark}

\newcommand{\sectiontriviale}{s}
\newcommand{\lineaire}{\Theta}
\def\build#1_#2^#3{\mathrel{\mathop{\kern0pt#1}\limits_{#2}^{#3}}}

\begin{document}
\begin{abstract}
Let $M$ be any compact simply-connected $d$-dimensional smooth manifold
and let $\mathbb{F}$ be any field.
We show that the Gerstenhaber algebra structure on the Hochschild
cohomology on the singular cochains of $M$, $HH^*(S^*(M);S^*(M))$,
extends to a Batalin-Vilkovisky algebra.
Such Batalin-Vilkovisky algebra was conjecturated to exist
and is expected to be isomorphic to the Batalin-Vilkovisky
algebra on the free loop space homology on $M$, $H_{*+d}(LM)$ introduced
by Chas and Sullivan.
We also show that the negative cyclic cohomology
$HC^*_-(S^*(M))$ has a Lie bracket.
Such Lie bracket is expected to coincide with the Chas-Sullivan string
bracket on the equivariant homology $H_*^{S^1}(LM)$.
\end{abstract}
\title{\bf Batalin-Vilkovisky algebra structures on Hochschild Cohomology.}
\author{Luc Menichi}
\address{UMR 6093 associ\'ee au CNRS\\
Universit\'e d'Angers, Facult\'e des Sciences\\
2 Boulevard Lavoisier\\49045 Angers, FRANCE}
\email{firstname.lastname at univ-angers.fr}
\keywords{String Topology, Batalin-Vilkovisky algebra, Gerstenhaber
  algebra, Hochschild cohomology, free loop space}
\maketitle
\begin{center}
\end{center}
\section{Introduction}
Except where specified, we work over an arbitrary field $\mathbb{F}$.
Let $M$ be a compact oriented $d$-dimensional smooth manifold.
Denote by $LM:=map(S^1,M)$ the free loop space on $M$.
Chas and Sullivan~\cite{Chas-Sullivan:stringtop} have shown that the shifted free loop homology
$H_{*+d}(LM)$ has a structure of Batalin-Vilkovisky algebra
(Definition~\ref{definition BV algebre}). In particular, they showed
that $H_{*+d}(LM)$ is a Gerstenhaber algebra (Definition~\ref{definition algebre de Gerstenhaber}).
On the other hand, let $A$ be a differential graded algebra.
The Hochschild cohomology of $A$ with coefficients in $A$, $HH^*(A;A)$,
is a Gerstenhaber algebra.
These two Gerstenhaber algebras are expected to be related:
\begin{conjecture}(due to ~\cite[``dictionary'' p. 5]{Chas-Sullivan:stringtop} or~\cite{CohJon:ahtrostringtopology}?)\label{conjecture iso algebre de Gerstenhaber}
\noindent If $M$ is simply connected then there is an isomorphism of Gerstenhaber algebras
$H_{*+d}(LM)\cong HH^*(S^*(M);S^*(M))$ between the free loop space homology and the Hochschild cohomology of the
algebra of singular cochains on $M$.
\end{conjecture}
F\'elix, Thomas and 
Vigu\'e-Poirrier~\cite[Appendix]{Felix-Thomas-Vigue:Hochschildmanifold}
proved that there is a linear isomorphism of lower degree $d$
\begin{equation}\label{iso FTV}
 \mathbb{D}:HH^{-p-d}(S^*(M),S^*(M)^\vee)\buildrel{\cong}\over
\rightarrow HH^{-p}(S^*(M),S^*(M)).
\end{equation}
We prove
\begin{theor}
(Theorem~\ref{structure BV sur la cohomologie de Hochschild de M})
The Connes coboundary map on 
$HH^*(S^*(M),S^*(M)^\vee)$ defines via the isomorphism~(\ref{iso FTV})
a structure of Batalin-Vilkovisky algebra on
the Gerstenhaber algebra $HH^*(S^*(M),S^*(M))$.
\end{theor}

Assume that $M$ is simply-connected.
Jones~\cite{JonesJ:Cycheh} proved that there is an isomorphism
$$J:H_{p+d}(LM)\buildrel{\cong}\over\rightarrow
HH^{-p-d}(S^*(M),S^*(M)^\vee)$$
such that the $\Delta$ operator of the Batalin-Vilkovisky algebra
$H_{*+d}(LM)$ and Connes coboundary map $B^\vee$ on
$HH^{*-d}(S^*(M),S^*(M)^\vee)$ satisfies $J\circ\Delta=B^\vee\circ J$.
Of course, we conjecture:
\begin{conjecture}\label{conjecture iso d'algebres}
The isomorphism
$$\mathbb{D}\circ J:H_{p+d}(LM)\buildrel{\cong}\over\rightarrow
 HH^{-p}(S^*(M),S^*(M))$$
is a morphism of graded algebras.
\end{conjecture}
Notice that Conjecture~\ref{conjecture iso d'algebres}
implies that the composite $\mathbb{D}\circ J$ is an isomorphism
of Batalin-Vilkovisky algebras between the Chas-Sullivan
Batalin-Vilkovisky algebra and the Batalin-Vilkovisky  algebra
defined by Theorem~\ref{structure BV sur la cohomologie de Hochschild de M}.
Therefore Conjecture~\ref{conjecture iso d'algebres}
implies Conjecture~\ref{conjecture iso algebre de Gerstenhaber}.

Cohen and Jones~\cite[Theorem 3]{CohJon:ahtrostringtopology}
have an isomorphism of algebras
$$H_{p+d}(LM)\buildrel{\cong}\over\rightarrow
 HH^{-p}(S^*(M),S^*(M)).$$
So one should check perhaps if the isomorphism of Cohen-Jones
coincides with the isomorphism $\mathbb{D}\circ J$.
Over the reals or over the rationals, two proofs of such
an isomorphism of graded algebras
have been given by Merkulov~\cite{Merkulov:DeRhammfst}
and F\'elix, Thomas, Vigu\'e-Poirrier~\cite{Felix-Vigue-Thomas:ratstringtop}.

Theorem~\ref{structure BV sur la cohomologie de Hochschild de M}
comes from a general result
(Propositions~\ref{iso cohomologie homologie de Hochschild} and
\ref{BV algebre cohomologie de Hochschild})
who shows that the Hochschild cohomology $HH^*(A;A)$ of a differential
graded algebra $A$ which is a ``homotopy symmetric algebra'',
is a Batalin-Vilkovisky algebra.
As second application of this general result,
we recover the following theorem due to Thomas Tradler.
\begin{theor}~\cite[Example 2.15 and Theorem 3.1]{Tradler:bvalgcohiiip}
(Corollary~\ref{Hochschild algebre symmetrique non graduee})
Let $A$ be a symmetric algebra.
Then $HH^{*}(A,A)$ is
a Batalin-Vilkovisky algebra.
\end{theor}
This theorem has been reproved and extended by many
people~\cite{MenichiL:BValgaccoHa,Kaufmann:procvDelignecvc,Tradler-Zeinalian:ontcyclicDelignec,Costello:tcftcalabiyaucat,Kaufmann:modsahcfaco,Kontsevich-Soibelman:noainfinityalg,PoHu:HochschildcohologyPoincare,Eu-Schedler:CYFalgebra}
(in chronological order).
The last proof, the proof of Eu et Schedler~\cite{Eu-Schedler:CYFalgebra}
looks similar to ours.

Thomas Tradler gave a somehow complicated proof of the previous theorem
(Corollary~\ref{Hochschild algebre symmetrique non graduee}).
Indeed, his goal was to prove our main theorem
(Theorem~\ref{structure BV sur la cohomologie de Hochschild de M}).
In~\cite{Tradler-Zeinalian:InfinityPoincareduality}
or in~\cite{Tradler-Zeinalian:algebraicstringoperations},
Tradler and Zeinalian proved
Theorem~\ref{structure BV sur la cohomologie de Hochschild de M}
but only over a field of characteristic $0$
~\cite[``rational simplicial chain'' in the abstract]{Tradler-Zeinalian:InfinityPoincareduality} or~\cite[Beginning of 3.1]{Tradler-Zeinalian:algebraicstringoperations}.
Costello's result~\cite[Section 2.1]{Costello:tcftcalabiyaucat}
is also over a field of characteristic $0$.

Over $\mathbb{Q}$, we explain how to put
a Batalin-Vilkovisky algebra structure on
$
HH^*(S^*(M;\mathbb{Q});S^*(M;\mathbb{Q}))
$
(Corollary ~\ref{Hochschild rationnellement}) 
from a slight generalisation of
Corollary~\ref{Hochschild algebre symmetrique non graduee}
(Theorem~\ref{Hochschild algebre symmetrique}).
In fact both F\'elix, Thomas~\cite{Felix-Thomas:ratBVstringtop} and
Chen~\cite[Theorem 5.4]{Chen:generalchainmodel} proved
that the Chas-Sullivan
Batalin-Vilkovisky algebra $H_{p+d}(LM;\mathbb{Q})$
is isomorphic to the Batalin-Vilkovisky algebra
given by Corollary~\ref{Hochschild rationnellement}.

Finally, remark, that, over $\mathbb{Q}$, when the manifold $M$
is formal, a consequence of  F\'elix and Thomas work
~\cite{Felix-Thomas:ratBVstringtop}, is that 
$H_{p+d}(LM)$
is always isomorphic to the Batalin-Vilkovisky algebra
$HH^*(H^*(M);H^*(M)^\vee)$
given by Corollary~\ref{Hochschild algebre symmetrique non graduee}
applied to symmetric algebra $H^*(M)$.
Over $\mathbb{F}_2$, in~\cite{menichi:stringtopspheres},
we showed that this is not the case.
The present paper seems to explain why:

The Batalin-Vilkovisky algebra on $HH^*(S^*(M);S^*(M))$
given by Theorem~\ref{structure BV sur la cohomologie de Hochschild de M}
depends of course of the algebra $S^*(M)$ but also of a fundamental
class $[m]\in HH^*(S^*(M);S^*(M)^\vee)$ which seems hard to compute.
This fundamental class $[m]$ involves chain homotopies
for the commutativity of the algebra  $S^*(M)$.

The Batalin-Vilkovisky algebra on
$HH^*(S^*(M;\mathbb{Q});S^*(M;\mathbb{Q}))$
given by Corollary~\ref{Hochschild rationnellement}, depends of

-a commutative algebra,
Sullivan's cochain algebra of polynomial differential forms $A_{PL}(M)$
~\cite{Felix-Halperin-Thomas:ratht},

-and of the fundamental class $[M]\in H(A_{PL}(M))$.

{\em Acknowledgment:}
We wish to thank Jean-Claude Thomas for a discussion concerning
Ginzburg's preprint~\cite{Ginzburg:CYalgebras}.
We would like also to thank Yves F\'elix for explaining us the wonderful
isomorphism~(\ref{iso FTV}).
\section{Hochschild homology and cohomology}

Let $A$ be a differential graded algebra. 
Denote by $sA$ the suspension of $A$, $(s A)_i=A_{i-1}$.
Let $d_1$ be the differential on the tensor product of complexes
$A\otimes T(sA)\otimes A$.
We denote the tensor product of the elements $a\in A$, $sa_1\in sA$,
\ldots , $sa_k\in sA$ and $b\in A$ by $a[a_1|\cdots|a_k]b$.
Let $d_2$ be the differential on the graded vector space
$A\otimes T(sA)\otimes A$ defined by:
\begin{align*}
d_2a[a_1|\cdots|a_k]b=&(-1)^{\vert a\vert } aa_1[a_2|\cdots|a_k]b\\
&+\sum _{i=1}^{k-1} (-1)^{\varepsilon_i}{a[a_1|\cdots|a_ia_{i+1}|\cdots|a_k]b}\\
&-(-1)^{\varepsilon_{k-1}} a[a_1|\cdots|a_{k-1}]a_kb;
\end{align*}
Here $\varepsilon_i=\vert a\vert +\vert a_1\vert+\cdots +\vert a_i\vert+i$.
The {\it bar resolution of $A$}, denoted
$B(A;A;A)$, is the differential graded $(A,A)$-bimodule $(A\otimes
T(sA)\otimes A,d_1+d_2)$.

Denote by $A^{op}$ the opposite algebra of $A$.
Recall that any $(A,A)$-bimodule can be considered as a left (or right)
$A\otimes A^{op}$-module.
The {\it Hochschild chain complex} is the complex
$A\otimes_{A\otimes A^{op}} B(A;A;A)$ denoted $\mathcal{C}_*(A;A)$.
Explicitly $\mathcal{C}_*(A;A)$ is the complex ($A\otimes T(sA),d_1+d_2)$
with $d_1$ obtained by tensorization and
\begin{align*}
d_2a[a_1|\cdots|a_k]=&(-1)^{\vert a\vert } aa_1[a_2|\cdots|a_k]\\
&+\sum _{i=1}^{k-1} (-1)^{\varepsilon_i}{a[a_1|\cdots|a_ia_{i+1}|\cdots|a_k]}\\
&-(-1)^{\vert sa_k\vert\varepsilon_{k-1}} a_ka[a_1|\cdots|a_{k-1}].
\end{align*}
The {\it Hochschild homology} is the homology of the Hochschild chain complex:
$$HH_*(A;A):=H(\mathcal{C}(A;A)).$$
Let $M$ be a differential graded $(A,A)$-bimodule.
 The {\it Hochschild cochain complex} of $A$ with coefficients in $M$ is
the   complex
    $$
    \mathcal{C}^\ast  (A;M) =  (\mbox{Hom}   (T(sA), M), D_0+D_1) \,.
    $$
Here for $f \in  \mbox{Hom}( T(sA), M) $, 
  $D_0(f)([\,]) = d_M(f([\,]))$, $D_1(f)([\,])=0$, and for $k\geq 1$,
we have:
   $$D_0(f)([a_1|a_2|...|a_k])
  = d_{M}(f\left([a_1|a_2|...|a_k])\right) - \sum _{i=1} ^k
  (-1)^{\overline \epsilon_i} f([a_1|...|d_Aa_i|...|a_k])$$  and
  $$\renewcommand{\arraystretch}{1.6}
\begin{array}{ll}
D_1(f)([a_1|a_2|...|a_k])= &
   - (-1)^{|sa_1|\, |f|}  a_1 f([a_2|...|a_k])\\ &
- \sum _{i=2} ^k
  (-1)^{\overline \epsilon_i} f([a_1|...|a_{i-1}a_i|...|a_k])\\
 & + (-1) ^{\overline\epsilon _k} f([a_1|a_2|...|a_{k-1}])a_k
  \,,
\end{array}
\renewcommand{\arraystretch}{1}
$$
 where
  $\overline \epsilon _i = |f|+|sa_1|+|sa_2|+...+|sa_{i-1}|$.

The {\it Hochschild cohomology of $A$ with coefficients in $M$} is
 $$
 HH^*(A;M)=  H( \mathcal{C} ^\ast(A;M)) =H(\mbox{Hom}   (T(A), M), D_0+D_1)\,.
 $$
Since we work over an arbitrary field $\mathbb{F}$,
the bar resolution $B(A;A;A)\buildrel{\simeq}\over\rightarrow A$ is a
semi-free resolution of $A$ as an
$(A,A)$-bimodule~\cite[Proposition 19.2(ii)]{Felix-Halperin-Thomas:ratht}.
Therefore the Hochschild homology of $A$ is the differential torsion product
$$
HH_*(A;A)=\mbox{Tor}^{A\otimes A^{op}}(A,A)
$$
and the Hochschild cohomology is
 $$
 HH^*(A;M)\cong H(\mbox{Hom}_{A\otimes A^{op}}(B(A;A;A),M))=\mbox{Ext}_{A\otimes A^{op}}(A,M)
$$
where the latter denotes the differential "Ext" in the sense of
J.C. Moore (cf~\cite[Appendix]{Felix-Halperin-Thomas:Gors}).

Gerstenhaber proved that the Hochschild cohomology of $A$
with coefficients in $A$, $HH^*(A;A)$, is a Gerstenhaber
algebra~\cite{Gerstenhaber:cohosar}.
\begin{defin}\label{definition algebre de Gerstenhaber}

A {\it Gerstenhaber algebra} is a
commutative graded algebra $A$
equipped with a linear map
$\{-,-\}:A_i \otimes A_j \to A_{i+j+1}$ of degree $1$
such that:

\noindent a) the bracket $\{-,-\}$ gives $A$ a structure of graded
Lie algebra of degree $1$. This means that for each $a$, $b$ and $c\in A$

$\{a,b\}=-(-1)^{(\vert a\vert+1)(\vert b\vert+1)}\{b,a\}$ and 

$\{a,\{b,c\}\}=\{\{a,b\},c\}+(-1)^{(\vert a\vert+1)(\vert b\vert+1)}
\{b,\{a,c\}\}.$

\noindent b)  the product and the Lie bracket satisfy the following relation
called the Poisson relation:
$$\{a,bc\}=\{a,b\}c+(-1)^{(\vert a\vert+1)\vert b\vert}b\{a,c\}.$$
\end{defin}
In this paper, we show that for some algebras $A$, the Gerstenhaber
algebra structure of $HH^*(A;A)$ extends to a Batalin-Vilkovisky algebra.
\begin{defin}\label{definition BV algebre}
A {\it Batalin-Vilkovisky algebra} is a Gerstenhaber algebra
$A$ equipped with a degree $1$ linear map $\Delta:A_{i}\rightarrow A_{i+1}$
such that $\Delta\circ\Delta=0$ and
\begin{equation}\label{relation BV} 
\{a,b\}=(-1)^{\vert a\vert}\left(\Delta(a\cup b)-(\Delta a)\cup b-(-1)^{\vert
  a\vert}a\cup(\Delta b)\right)
\end{equation}
for $a$ and $b\in A$.
\end{defin}
\section{The isomorphism between $HH^*(A;A)$ and $HH^*(A;A^\vee)$}
In this section, we first present a method that gives an isomorphism
between the Hochschild cohomology of $A$ with coefficients in $A$,
$HH^*(A;A)$ and the Hochschild cohomology of $A$ with coefficients in
the dual $A^\vee$, $HH^*(A;A^\vee)$.
This method is a generalisation of the method used by F\'elix, Thomas
and Vigu\'e-Poirrier to obtain the isomorphism~(\ref{iso FTV}).
Then we show that this isomorphism looks like a Poincar\'e duality
isomorphism: this isomorphism is given by the action of the algebra
$HH^*(A;A)$ on a fundamental class $[m]\in HH^*(A;A^\vee)$. 

Let us first recall the definition of the action
of $HH^*(A;A)$ on $HH^*(A;A^\vee)$.
Let $A$ be a (differential graded) algebra.
Let $M$ and $N$ two $A$-bimodules.
Let $f\in\mathcal{C}^*(A,M)$ and $g\in\mathcal{C}^*(A,N)$.
We denote by $\otimes_A(f,g)\in\mathcal{C}^*(A,M\otimes_A N)$
the linear map defined by
$$\otimes_A(f,g)([a_1|\dots|a_n])=
\sum_{p=0}^n f([a_1|\dots|a_p])\otimes_A g([a_{p+1}|\dots|a_n]).
$$
This define a natural morphism of complexes
$$
\otimes_A:\mathcal{C}^*(A,M)\otimes \mathcal{C}^*(A,N)
\rightarrow \mathcal{C}^*(A,M\otimes_A N)
$$
Therefore, in homology, we have a natural morphism
$$
H_*(\otimes_A):HH^*(A,M)\otimes HH^*(A,N)
\rightarrow HH^*(A,M\otimes_A N)
$$
If we let take $A=M$,
and use the isomorphism of $A$-bimodules
$$A\otimes_A N\buildrel{\cong}\over\rightarrow N, a\otimes_A n\mapsto a.n,$$
the composite
\begin{equation}\label{action en cohomologie de Hochschild}
\mathcal{C}^*(A,A)\otimes \mathcal{C}^*(A,N)
\buildrel{\otimes_A}\over\rightarrow \mathcal{C}^*(A,A\otimes_A N)\cong \mathcal{C}^*(A,N)
\end{equation}
is a left action of $\mathcal{C}^*(A,A)$ on $\mathcal{C}^*(A,N)$.
In the particular case, $A=M=N$, this composite is the usual
cup product on $\mathcal{C}^*(A,A)$ denoted $\cup$.

Denote by $A^\vee$ the dual of $A$.
Let $\eta:\mathbb{F}\rightarrow A$ be the unit of the algebra.
Then we have a natural map 
$$HH^*(\eta, A^\vee):HH^*(A,A^\vee)\rightarrow
HH^*(\mathbb{F},A^\vee)\cong H(A^\vee).$$
\begin{proposition}\label{iso cohomologie homologie de Hochschild}
Let $[m]\in HH^{-d}(A,A^\vee)$ be an element of lower degree $d$ such that
the morphism of left $H(A)$-modules
$$H(A)\buildrel{\cong}\over
\rightarrow H(A^\vee), a\mapsto a.HH^{-d}(\eta, A^\vee)([m])$$
is an isomorphism.
Then the action of $HH^*(A,A)$ on $[m]\in HH^{-d}(A,A^\vee)$
gives the isomorphism of lower degree $d$ of $HH^*(A,A)$-modules
$$HH^p(A,A)\buildrel{\cong}\over\rightarrow
HH^{p-d}(A,A^\vee),\;a\mapsto a\cdot[m].$$
\end{proposition}
\begin{proof}
Let $\varepsilon_A:P\buildrel{\simeq}\over\rightarrow A$
be a resolution of $A$ as left $A\otimes A^{op}$-semifree module.
Let $s_A:A\buildrel{\simeq}\over\hookrightarrow P$ be a morphism of
left $A$-modules which is a
section of $\varepsilon_A$.
The morphism $HH^*(\eta, A^\vee):HH^*(A,A^\vee)\rightarrow HH^*(\mathbb{F},A^\vee)$
is equal to the following composite of
$$
HH^*(A,A^\vee):=Ext_{A\otimes A^{op}}(A,A^\vee)
\buildrel{Ext_{i_1}(A,A^\vee)}\over\rightarrow
Ext_{A}(A,A^\vee)
$$
and
$$
Ext_{A}(A,A^\vee)
\build\rightarrow_\cong^{Ext_{\eta}(\eta,A^\vee)}
Ext_\mathbb{F}(\mathbb{F},A^\vee)=:
HH^*(\mathbb{F},A^\vee)
$$
where $i_1:A\hookrightarrow A\otimes A^{op}$
is the inclusion of the first factor.

Therefore, $HH^*(\eta, A^\vee)$ is the map induced in homology
by the composite
$$
Hom_{A\otimes A^{op}} (P,A^\vee)\buildrel{Hom(s_A,A^\vee)}
\over\rightarrow Hom_{A} (A,A^\vee)\build\rightarrow_\cong^{ev(1_A)} A^\vee.
$$
where $ev(1_A)$ is the evaluation at the unit $1_A\in A$.
This composite maps the cycle $m\in Hom_{A\otimes A^{op}} (P,A^\vee)$
to $m\circ s_A$ and then to $(m\circ s_A)(1_A)$.
Since $m\circ s_A:A\rightarrow A^\vee$ maps $a\in A$ to
$a\cdot \left((m\circ s_A)(1)\right)$,
by hypothesis, $m\circ s_A$ is a quasi-isomorphism.
Since $s_A$ is a quasi-isomorphism,
$m:P\buildrel{\simeq}\over\rightarrow A^{\vee}$ is also a quasi-isomorphism.

By applying the functor $Hom_{A\otimes A^{op}} (P,-)$ to
the two quasi-isomorphisms of $A$-bimodules
$$A\build\leftarrow_{\simeq}^{\varepsilon_A}
P\build\rightarrow_{\simeq}^{m}
A^\vee,$$
we obtain the quasi-isomorphism of complexes
$$Hom_{A\otimes A^{op}} (P,A)
\build\leftarrow_{\simeq}^{\varepsilon_A}
Hom_{A\otimes A^{op}} (P,P)
\build\rightarrow_{\simeq}^{Hom_{A\otimes A^{op}} (P,m)}
Hom_{A\otimes A^{op}} (P,A^\vee).
$$
By applying homology, we get the desired isomorphism,
since the action of $HH^*(A,A)$ on $HH^*(A,A^\vee)$
is induced by the composition map
$$
Hom_{A\otimes A^{op}} (P,A^\vee)\otimes Hom_{A\otimes A^{op}} (P,P)
\rightarrow Hom_{A\otimes A^{op}} (P,A^\vee)
$$
$$
m\otimes f\mapsto m\circ f=Hom_{A\otimes A^{op}} (P,m)(f)
$$
Alternatively, the two isomorphisms
$$HH^*(A,A)\build\leftarrow_{\cong}^{HH*(A,\varepsilon_A)}
HH^*(A,P)\build\rightarrow_{\cong}^{HH^*(A,m)} HH^*(A,A^\vee)$$
maps $\varepsilon_A$ (which is the unit of $HH^*(A,A)$) to $id_P:P\rightarrow P$ and then to $m$.
They are morphisms of $HH^*(A,A)$-modules
since 
$$H_*(\otimes_A):HH^*(A,A)\otimes HH^*(A,N)
\rightarrow HH^*(A,A\otimes_A N)$$ is natural with respect to $N$.
\end{proof}
\section{Batalin-Vilkovisky algebra structures on Hochschild cohomology}
In this section, we explain when an isomorphism
$HH^*(A;A)\cong HH^*(A;A^\vee)$ gives a Batalin-Vilkovisky algebra structure
on the Gerstenhaber algebra $HH^*(A;A)$.
Our proof relies on the proof of a similar result due to
Ginzburg~\cite[Theorem 3.4.3 (ii)]{Ginzburg:CYalgebras}.
Ginzburg basically explains when an isomorphism
$HH^*(A;A)\cong HH_*(A;A)$ gives a Batalin-Vilkovisky algebra
structure on $HH^*(A;A)$.

Denote by $B$ Connes boundary in the Hochschild complex $\mathcal{C}_*(A;A)$
and by $B^\vee$ its dual in
$\mathcal{C}^*(A,A^\vee)\cong\mathcal{C}_*(A;A)^\vee$. We prove:
\begin{proposition}\label{BV algebre cohomologie de Hochschild}
Let $[m]\in HH^{-d}(A,A^\vee)$ such that
the morphism of  $HH^*(A,A)$-modules
$$HH^p(A,A)\buildrel{\cong}\over\rightarrow
HH^{p-d}(A,A^\vee),\;a\mapsto a.[m]$$
is an isomorphism.
If $H_*(B^\vee)([m])=0$ then the Gerstenhaber algebra
$HH^*(A,A)$ equipped with $H_*(B^\vee)$ is a Batalin-Vilkovisky algebra.
\end{proposition}
As we will see Proposition~\ref{BV algebre cohomologie de Hochschild}
is almost the dual of the following Proposition due to Victor Ginzburg.
Recall first that the Hochschild cohomology of a (differential graded)
algebra, acts on its Hochschild homology

$$HH^p(A;A)\otimes HH_d(A;A)\rightarrow HH_{p-d}(A;A)
$$
$$\eta\otimes c\mapsto i_\eta(c)=\eta. c$$
In non-commutative geometry, the action of $\eta\in HH^*(A;A)$
on $c\in HH_*(A;A)$ is denoted by $i_{\eta}(c)$.
\begin{proposition}~\cite[Theorem 3.4.3 (ii)]{Ginzburg:CYalgebras}\label{BV algebre homologie de Hochschild}
Let $c\in HH_{d}(A,A)$ such that
the morphism of  $HH^*(A,A)$-modules
$$HH^p(A,A)\buildrel{\cong}\over\rightarrow
HH_{d-p}(A,A),\;\eta\mapsto\eta.c$$
is an isomorphism.
If $H_*(B)(c)=0$ then the Gerstenhaber algebra
$HH^*(A,A)$ equipped with $H_*(B)$ is a Batalin-Vilkovisky algebra.
\end{proposition}
\begin{rem}
The condition $H_*(B)(c)=0$ does not appear
in~\cite[Theorem 3.4.3 (ii)]{Ginzburg:CYalgebras} since according to
Ginzburg, this condition is automatically satisfied for a Calabi-Yau
algebra of dimension $d$.
In both Propositions~\ref{BV algebre cohomologie de Hochschild} and
\ref{BV algebre homologie de Hochschild}, if the condition
 $H_*(B^\vee)([m])=0$ or  $H_*(B)(c)=0$ is not satisfied,
$\Delta(1)$ can be non zero and the relation~(\ref{relation BV}) is replaced
by the more general relation
\begin{multline*}
\{\xi,\eta\}=(-1)^{\vert \xi\vert} [\Delta(\xi\cup\eta)-\\
(-1)^{\vert \xi\vert} \xi\cup(\Delta\eta)
-(\Delta\xi)\cup\eta+(-1)^{\vert \xi\vert +\vert \eta\vert }\xi\cup\eta\cup(\Delta 1)].
\end{multline*}
\end{rem}
\begin{proof}[Proof of Proposition~\ref{BV algebre homologie de Hochschild}]
By definition the $\Delta$ operator on $HH^*(A;A)$ is given
by $(\Delta a).c:=B(a.c)$ for any $a\in HH^*(A;A)$.
Therefore the proposition follows from the following Lemma
due to Victor Ginzburg.
\end{proof}
\begin{lem}~\cite[formula (9.3.2)]{Ginzburg:CYalgebras}\label{lemme homologie de Hoschchild}
Let $A$ be a differential graded algebra.
For any $\eta,\xi\in HH^*(A;A)$ and $c\in HH_*(A;A)$,
\begin{multline*}
\{\xi,\eta\}.c=(-1)^{\vert \xi\vert}
B\left[(\xi\cup\eta).c\right]-\xi.B(\eta.c)\\
+(-1)^{(\vert \eta\vert +1)(\vert \xi\vert +1)}\eta.B(\xi.c)
+(-1)^{\vert \eta\vert} (\xi\cup\eta).B(c).
\end{multline*}
\end{lem}
\begin{proof}
Let us recall the proof of Victor Ginzburg.
Denote by
$$HH^p(A;A)\otimes HH_j(A;A)\rightarrow HH_{j-p+1}(A;A)$$
$$(\eta,a)\mapsto L_{\eta}(a)
$$
the action of the suspended graded Lie algebra $HH^*(A;A)$
on $HH_*(A;A)$.
Gelfand, Daletski and
Tsygan~\cite{Daletski-Gelfand-Tsygan:variantnoncommdiffgeom} proved that the Gerstenhaber algebra
$HH^*(A;A)$ and Connes boundary map $B$ on $HH_*(A;A)$ form a calculus
~\cite[p. 93]{Cuntz-Skandalis-Tamarkin:cyclichomnoncomgeom}.
Therefore, we have the following equalities
\begin{multline*}
i_{\{\xi,\eta\}}=\{L_\xi,i_\eta\}=
L_\xi\circ i_\eta-(-1)^{(\vert \xi\vert +1)\vert \eta\vert }
i_{\eta}\circ L_\xi\\
=(-1)^{\vert \xi\vert }\{B,i_{\xi}\}\circ i_\eta
-(-1)^{(\vert \xi\vert +1)\vert \eta\vert }i_{\eta}
\circ(-1)^{\vert \xi\vert} \{B,i_{\xi}\}\\
=(-1)^{\vert \xi\vert } B\circ i_{\xi}\circ i_\eta
-i_{\xi}\circ B\circ i_\eta
+(-1)^{(\vert \eta\vert +1)(\vert \xi\vert +1)}i_{\eta}\circ B\circ i_{\xi}
+(-1)^{\vert \eta\vert (\vert \xi\vert +1)}i_{\eta}\circ i_{\xi}\circ B\\
=(-1)^{\vert \xi\vert} B\circ i_{\xi\cup\eta}-i_{\xi}\circ B\circ i_{\eta}
+(-1)^{(\vert\eta\vert+1)(\vert\xi\vert+1)}i_{\eta}\circ B\circ i_{\xi}
+(-1)^{\vert\eta\vert} i_{\xi\cup\eta}\circ B.
\end{multline*}
By applying this equality of operators to $c$, we obtain the Lemma.
\end{proof}
We now prove the following Lemma which is the dual of Lemma~\ref{lemme homologie de Hoschchild}.
\begin{lem}\label{lemme cohomologie de Hoschchild}
Let $A$ be a differential graded algebra.
For any $\eta,\xi\in HH^*(A;A)$ and $m\in HH^*(A;A^\vee)$,
\begin{multline*}
\{\xi,\eta\}.m=(-1)^{\vert\xi\vert}
B^\vee\left[(\xi\cup\eta).m\right]-\xi.B^\vee(\eta.m)\\
+(-1)^{(\vert\eta\vert+1)(\vert\xi\vert+1)}\eta.B^\vee(\xi.m)
+(-1)^{\vert\eta\vert}(\xi\cup\eta).B^\vee(m).
\end{multline*}
\end{lem}
\begin{proof}
The action of $HH^*(A;A)$ on $HH_*(A;A)$ comes from a (right) action
of the $\mathcal{C}^*(A;A)$ on $\mathcal{C}_*(A;A)$ given by
$$\mathcal{C}_*(A;A)\otimes\mathcal{C}^*(A;A)\rightarrow \mathcal{C}_*(A;A)$$
$$(m[a_1|\dots|a_n],f)\mapsto
i_f(m[a_1|\dots|a_n]):=\sum_{p=0}^n (m.f[a_1|\dots|a_p])[a_{p+1}|\dots|a_n].$$
Therefore $\mathcal{C}^*(A;A)$ acts on the left on the dual
$\mathcal{C}_*(A;A)^\vee$. Explicitly, the action is given by
$$\mathcal{C}^*(A;A)\otimes\mathcal{C}_*(A;A)^\vee
\rightarrow \mathcal{C}_*(A;A)^\vee$$
$$(f,\varphi)\mapsto \varphi\circ i_f.$$
Through the canonical isomorphism
$\mathcal{C}(A;A^\vee)\buildrel{\cong}\over\rightarrow\mathcal{C}_*(A;A)^\vee$,
$g\mapsto \varphi$ defined by
$\varphi(m[a_1|\dots|a_n]):=(g[a_1|\dots|a_n])(m)$, this left action coincides
with the left action defined by the
composite~(\ref{action en cohomologie de Hochschild}).

Let us precise our sign convention: we define $B^\vee$ by
$B^\vee(m):=(-1)^{\vert m\vert} m\circ B$.
Denote by $\varepsilon$ the sign
$(-1)^{\vert m\vert (\vert \xi\vert +\vert \eta\vert +1)}$.
For any $m\in HH_*(A;A)^\vee$, we have the following equalities:
$$
m(\{\xi,\eta\}.c)=\varepsilon(\{\xi,\eta\}.m)(c),
$$
$$
(-1)^{\vert \xi\vert } m\circ B[(\xi\cup\eta).c]=
(-1)^{\vert \xi\vert +\vert m\vert }[B^\vee(m)][(\xi\cup\eta).c]=
\varepsilon (-1)^{\vert \eta\vert }[(\xi\cup\eta).B^\vee(m)](c),
$$
\begin{multline*}
-m[\xi.B(\eta.c)]=
(-1)^{1+\vert m\vert \vert \xi\vert }[\xi.m]\circ B(\eta.c)=\\
(-1)^{1+\vert m\vert \vert \xi\vert +\vert \xi\vert +\vert m\vert }
[B^\vee(\xi.m)](\eta.c)=
\varepsilon (-1)^{(\vert \eta\vert +1)(\vert \xi\vert +1)}[\eta.B^\vee(\xi.m)](c),
\end{multline*}
by exchanging $\xi$ and $\eta$,
$$
(-1)^{(\vert \eta\vert +1)(\vert \xi\vert +1)} m[\eta.B(\xi.c)]=
-\varepsilon[\xi.B^\vee(\eta.m)](c),
$$
\begin{multline*}
(-1)^{\vert \eta\vert } m[(\xi\cup\eta).B(c)]=
\varepsilon(-1)^{\vert \eta\vert +\vert m\vert }[(\xi\cup\eta).m]\circ B(c)=\\
\varepsilon(-1)^{\vert \xi\vert} B^\vee[(\xi\cup\eta).m](c).
\end{multline*}
Therefore by evaluating the linear form m
on the terms of the equation
given by Lemma~\ref{lemme homologie de Hoschchild},
we obtain the desired equality.
\end{proof}
\begin{rem}
The equality in Lemma~\ref{lemme cohomologie de Hoschchild}
is the same as the equality in Lemma~\ref{lemme homologie de Hoschchild}.
In fact, alternatively, to prove Lemma~\ref{lemme cohomologie de Hoschchild},
we could have proved that the Gerstenhaber algebra
$HH^*(A;A)$ and the dual of Connes boundary map $B^\vee$ on $HH^*(A;A^\vee)$
form a calculus. Indeed, in the proof of
Lemma~\ref{lemme homologie de Hoschchild}, we have remarked that the
desired equality holds for any calculus.
\end{rem}
\begin{proof}[Proof of Proposition~\ref{BV algebre cohomologie de Hochschild}]
By definition the $\Delta$ operator on $HH^*(A;A)$ is given
by $(\Delta a).m:=B^\vee(a.m)$ for any $a\in HH^*(A;A)$.
Therefore the proposition follows from
Lemma~\ref{lemme cohomologie de Hoschchild}.
\end{proof}
\section{Applications}
As first application of
Proposition~\ref{BV algebre cohomologie de Hochschild},
we show
\begin{theor}\label{Hochschild algebre symmetrique}
Let $A$ be an algebra equipped with a degree $d$ quasi-isomorphism of
$A$-bimodules $\lineaire:A\buildrel{\simeq}\over\rightarrow A^{\vee}$
between $A$ and its dual $\text{Hom}(A,\mathbb{F})$.
Then the Connes coboundary map on $HH^{*}(A,A^{\vee})$ defines via
the isomorphism $HH^{*}(A,\lineaire):HH^{p}(A,A)
\buildrel{\cong}\over\rightarrow HH^{p-d}(A,A^{\vee})$
a structure of Batalin-Vilkovisky algebra on the Gerstenhaber algebra
$HH^{*}(A,A)$.
\end{theor}
In representation theory~\cite{Curtis-Reiner:represent}, an (ungraded) algebra
$A$ is {\it symmetric} if $A$ is equipped with an isomorphism of
$A$-bimodules $\lineaire:A\buildrel{\cong}\over\rightarrow A^{\vee}$
between $A$ and its dual $\text{Hom}(A,\mathbb{F})$.
The following Corollary is implicit in~\cite{Tradler:bvalgcohiiip}
and was for the first time explicited
in~\cite[Theorem 1.6]{MenichiL:BValgaccoHa}.
\begin{cor}~\cite{Tradler:bvalgcohiiip,
  MenichiL:BValgaccoHa}\label{Hochschild algebre symmetrique non graduee}
Let $A$ be a symmetric algebra.
Then $HH^{*}(A,A)$ is
a Batalin-Vilkovisky algebra.
\end{cor}
In \cite{Kaufmann:procvDelignecvc} or
\cite[Corollary 3.4]{Tradler-Zeinalian:ontcyclicDelignec} or
\cite[Section 1.4]{Costello:tcftcalabiyaucat} or
\cite[Theorem B]{Kaufmann:modsahcfaco} or
\cite[Section 11.6]{Kontsevich-Soibelman:noainfinityalg} or
\cite{PoHu:HochschildcohologyPoincare},
this Batalin-Vilkovisky algebra structure on $HH^{*}(A,A)$ extends
to a structure of algebra on the Hochschild cochain complex
$\mathcal{C}^{*}(A,A)$ over various operads or PROPs: the so-called
cyclic Deligne conjecture.
\begin{proof}[Proof of Theorem~\ref{Hochschild algebre symmetrique}]
Let $\varepsilon_A:P:=B(A;A;A)\buildrel{\simeq}\over\rightarrow A$ be
the bar resolution of $A$.
Denote by $m$ the composite
$P\build\rightarrow_\simeq^{\varepsilon_A}
A\build\rightarrow_\simeq^{\lineaire} A^\vee$.
Since $m$ commutes with the differential, $m$ is a cycle in
$Hom_{A\otimes A^{op}}(P,A)$.
As we saw in the proof of
Proposition~\ref{iso cohomologie homologie de Hochschild},
the composite $HH^*(A,m)\circ HH^*(A,\varepsilon_A)^{-1}:$
$$HH^*(A,A)\build\leftarrow_{\cong}^{HH*(A,\varepsilon_A)}
HH^*(A,P)\build\rightarrow_{\cong}^{HH^*(A,m)} HH^*(A,A^\vee)$$
coincides with the morphism of left $HH^*(A;A)$-modules
$$
HH^p(A;A)\buildrel{\cong}\over\rightarrow
HH^{p-d}(A;A^\vee),\; a\mapsto a\cdot m.
$$
By definition of $m$, this composite is also $HH^*(A,\lineaire)$.

Denote by by $\varepsilon_{BA}:TsA\twoheadrightarrow \mathbb{F}$
the canonical projection whose kernel is $T^+sA$.
Since $\varepsilon_A:B(A;A;A)\twoheadrightarrow A$ is the composite of
$A\otimes\varepsilon_{BA}\otimes A$ and of the multiplication on $A$
$$A\otimes TsA\otimes A\twoheadrightarrow A\otimes\mathbb{F} \otimes A
\cong A\otimes A\twoheadrightarrow A,$$ 
the canonical isomorphisms of complexes
$$
Hom_{A\otimes A^{op}}(B(A;A;A),A^\vee)\cong\mathcal{C}^*(A;A^\vee)
\cong\mathcal{C}_*(A;A)^\vee
$$
map $m$ to the linear form on $\mathcal{C}_*(A;A)$:
$$\lineaire(1)\otimes\varepsilon_{BA}:A\otimes TsA\twoheadrightarrow
\mathbb{F}\otimes\mathbb{F}\cong \mathbb{F}.$$
Connes (normalized or not)
boundary map $B:\mathcal{C}_*(A;A)\rightarrow \mathcal{C}_*(A;A)$
factorizes through $A\otimes T^+sA$.
So $B^\vee(\lineaire(1)\otimes\varepsilon_{BA})=
\pm(\lineaire(1)\otimes\varepsilon_{BA})\circ B=0$.
Therefore, we can apply
Proposition~\ref{BV algebre cohomologie de Hochschild}

Remark: In the case of
Corollary~\ref{Hochschild algebre symmetrique non graduee},
$m$ correspond to a trace $\lineaire(1)\in \mathcal{C}^0(A;A^\vee)$.
Since $H(B^\vee):HH^p(A;A^\vee)\rightarrow HH^{p-1}(A;A^\vee)$ decreases
(upper) degrees and $HH^p(A;A^\vee)=0$ for $p<0$, it is obvious that
$H(B^\vee)(\lineaire(1))=0$.
\end{proof}
Working, with rational coefficients, we easily obtain
\begin{cor}~\cite{Tradler-Zeinalian:InfinityPoincareduality}
\label{Hochschild rationnellement}
The Hochschild cohomology
$$
HH^{*}(S^*(M;\mathbb{Q});S^*(M;\mathbb{Q}))\cong 
HH^{*-d}(S^*(M;\mathbb{Q});S^*(M;\mathbb{Q}))
$$
is a Batalin-Vilkovisky algebra.
\end{cor}
Tradler and Zeinalian~\cite{Tradler-Zeinalian:InfinityPoincareduality}
give a proof of this result.
Here is a shorter proof, although we don't claim that we have obtained
the same Batalin-Vilkovisky algebra.
It should not be difficult to see that the Batalin-Vilkovisky algebra
given by Corollary~\ref{Hochschild rationnellement} coincides
with the Batalin-Vilkovisky algebra given by our main 
theorem (Theorem~\ref{structure BV sur la cohomologie de Hochschild de M})
in the case of the field $\mathbb{Q}$. Therefore, one could deduce
Corollary~\ref{Hochschild rationnellement} from
Theorem~\ref{structure BV sur la cohomologie de Hochschild de M}.
But it is much more simple to give a separate proof of
Corollary~\ref{Hochschild rationnellement}.
As we would like to emphasize in this paper, the rational case is much
more simple than the case of a field $\mathbb{F}$
of characteristic $p$ different from $0$.
\begin{proof}[Proof of Corollary~\ref{Hochschild rationnellement}]
Since we are working over $\mathbb{Q}$,
there exists quasi-isomorphisms of
algebras~\cite[Corollary 10.10]{Felix-Halperin-Thomas:ratht}
$S^*(M;\mathbb{Q})\buildrel{\simeq}\over\rightarrow
D(M) \buildrel{\simeq}\over\leftarrow A_{PL}(M)$
where $A_{PL}(M)$ is a commutative (differential graded) algebra.
Since the Gerstenhaber algebra structure on
Hochschild cohomology is preserves by quasi-isomorphism of
algebras~\cite[Theorem 3]{Felix-Menichi-Thomas:GerstduaiHochcoh},
we obtain an isomorphism of Gerstenhaber algebras
$$
HH^*(S^*(M;\mathbb{Q});S^*(M;\mathbb{Q}))\cong HH^*(A_{PL}(M);A_{PL}(M)).
$$
Since $H(A_{PL}(M))\cong H^*(M;Q)$, Poincar\'e duality
induces an quasi-isomorphism of $A_{PL}(M)$-modules,
and so of $A_{PL}(M)$-bimodules,
since the algebra $A_{PL}(M)$ is commutative:
$$
A_{PL}(M)\build\rightarrow_\cong^{\cap [M]} A_{PL}(M)^\vee.
$$
By applying Theorem~\ref{Hochschild algebre symmetrique},
we obtain that
$$HH^*(A_{PL}(M);A_{PL}(M))\cong HH^{*-d}(A_{PL}(M);A_{PL}(M)^\vee)$$
is a Batalin-Vilkovisky algebra.
\end{proof}
In~\cite{Felix-Thomas:ratBVstringtop} and
~\cite[Theorem 5.4]{Chen:generalchainmodel}, it is shown that the
Batalin-Vilkovisky algebra $H_{p+d}(LM;\mathbb{Q})$ of Chas
and Sullivan is isomorphic to the Batalin-Vilkovisky algebra
on 
$$HH^{-p}(S^*(M;\mathbb{Q});S^*(M;\mathbb{Q}))\cong
HH^{-p-d}(A_{PL}(M);A_{PL}(M)^\vee).$$
given by Corollary~\ref{Hochschild rationnellement}.

Recall the following theorem due to F\'elix, Thomas and
Vigu\'e-Poirrier.
\begin{theor}~\cite[Appendix]{Felix-Thomas-Vigue:Hochschildmanifold}
\label{Felix-Thomas-Vigue}
Let $M$ be a compact connected oriented $d$-dimensional smooth manifold.
Then there is an isomorphism of lower degree $d$
$$\mathbb{D}^{-1}:HH^p(S^*(M),S^*(M))\buildrel{\cong}\over\rightarrow
HH^{p-d}(S^*(M),S^*(M)^\vee).$$
\end{theor}
As second application of
Propositions~\ref{iso cohomologie homologie de Hochschild}
and \ref{BV algebre cohomologie de Hochschild}, we will recover
the isomorphism of F\'elix, Thomas and Vigu\'e-Poirrier and prove
our main theorem:  
\begin{theor}\label{structure BV sur la cohomologie de Hochschild de M}
Let $M$ be a compact connected oriented $d$-dimensional smooth manifold.
Let $[M]\in H_d(M)$ be its fundamental class.
Then 

1) For any $a\in HH^*(S^*(M),S^*(M))$, the image of $a$
by $\mathbb{D}^{-1}$ is given by the action
of $a$ on $(J\circ H_*(\sectiontriviale))([M])$:

$$\mathbb{D}^{-1}(a)= a\cdot (J\circ H_*(\sectiontriviale))([M]).$$

2) The Gerstenhaber algebra structure
on $HH^*(S^*(M);S^*(M))$ and Connes coboundary map $H(B^\vee))$
on $HH^*(S^*(M);S^*(M)^\vee)$ defines via the isomorphism
$\mathbb{D}^{-1}$ a structure of Batalin-Vilkovisky algebra.
\end{theor}
Here $s$ denotes $\sectiontriviale:M\hookrightarrow LM$ the
inclusion of the constant loops into $LM$.
Recall that $J:H_*(LM)\rightarrow HH^*(S^*(M),S^*(M)^\vee)$ is the morphism
introduced by Jones in~\cite{JonesJ:Cycheh}. If $M$ is supposed to be simply connected,
then $J$ is an isomorphism.
\begin{proof}[Proof of Theorem~\ref{Felix-Thomas-Vigue}
and of Theorem~\ref{structure BV sur la cohomologie de Hochschild de M}]
We first follow
basically~\cite[Appendix]{Felix-Thomas-Vigue:Hochschildmanifold}.
Denote by $ev:LM\twoheadrightarrow LM$, $l\mapsto l(0)$ the evaluation map.
The morphism $J$ of Jones fits into the commutative triangle.
$$
\xymatrix{
H_*(LM)\ar[rr]^J\ar[dr]_{H_*(ev)}
&& HH^*(S^*(M),S^*(M)^\vee)\ar[dl]^{HH^*(\eta,S^*(M)^\vee)}\\
& H_*(M)
}
$$
Since $\sectiontriviale$ is a section of the evaluation map $ev$,
$J\circ H_*(\sectiontriviale)$ is a section of $HH^*(\eta,S^*(M)^\vee)$.
Therefore $HH^*(\eta,S^*(M)^\vee)\circ J\circ H_*(\sectiontriviale)([M])=[M]$.

By Poincar\'e duality, the composite of the two morphisms of $H^*(M)$-module
$$
H^*(M)\buildrel{\cap [M]}\over\rightarrow H_*(M)\cong H(S^*(M)^\vee),
a\mapsto a\cap [M]\mapsto a.[M].
$$
is an isomorphism of lower degree $d$. Therefore by applying
Proposition~\ref{iso cohomologie homologie de Hochschild} to
$[m]:=J\circ H_*(\sectiontriviale))([M]))$, we obtain
Theorem~\ref{Felix-Thomas-Vigue} and part 1) of
Theorem~\ref{structure BV sur la cohomologie de Hochschild de M}.

Consider $M$ equipped with the trivial $S^1$-action.
The section 
$\sectiontriviale:M\hookrightarrow LM$ is $S^1$-equivariant.
Therefore $\Delta\left(H_*(\sectiontriviale)([M])\right)=0$.
Recall that the Jones morphism $J$ satisfies $J\circ\Delta=H_*(B^\vee)\circ J$.
Therefore, since $(H_*(B^\vee)\circ J\circ H_*(\sectiontriviale))([M])=0$,
by applying Proposition~\ref{BV algebre cohomologie de Hochschild}, we obtain part 2) of Theorem~\ref{structure BV sur la cohomologie de Hochschild de M}.
\end{proof}
\begin{rem}
Part 1) of
Theorem~\ref{structure BV sur la cohomologie de Hochschild de M}
means exactly that
the morphism
$$\mathbb{D}^{-1}:HH^p(S^*(M),S^*(M))\buildrel{\cong}\over\rightarrow
HH^{p-d}(S^*(M),S^*(M)^\vee)$$ is the unique morphism
of $HH^*(S^*(M),S^*(M))$-modules such that 
the composite
$$J^{-1}\circ\mathbb{D}^{-1}:HH^{-p}(S^*(M),S^*(M))
\buildrel{\cong}\over\rightarrow H_{p+d}(LM)$$
respects the units of the algebras.
We conjecture (Conjecture~\ref{conjecture iso d'algebres})
that $J^{-1}\circ\mathbb{D}^{-1}$ respects also the products. 
\end{rem}
\section{cyclic homology}
In this section, we prove
\begin{cor}\label{crochet en cyclique}
Let $M$ be a compact oriented smooth $d$-dimensional manifold.
Then the negative cyclic cohomology on the singular cochains of $M$,
$HC^*_-(S^*(M))$, is a graded Lie algebra of lower degree $2-d$.
\end{cor}
If $M$ is simply-connected, Jones~\cite{JonesJ:Cycheh} proved
that there is an isomorphism
$$
H^{S^1}_*(LM)\buildrel{\cong}\over\rightarrow HC^*_-(S^*(M)).
$$
In~\cite{Chas-Sullivan:stringtop}, Chas and Sullivan defined
a Lie bracket, called the string bracket
$$
\{\; ,\; \}:
H^{S^1}_p(LM)\otimes H^{S^1}_q(LM)\rightarrow H^{S^1}_{p+q+2-d}(LM)
$$
Of course, we expect the two a priori different brackets to
be related:
\begin{conjecture}\label{conjecture iso d'algebres de Lie}
The Jones isomorphism $$
H^{S^1}_*(LM)\buildrel{\cong}\over\rightarrow HC^*_-(S^*(M)) 
$$ is an isomorphism of graded Lie algebras between Chas-Sullivan string
bracket and the Lie bracket defined in Corollary~\ref{crochet en cyclique}.
\end{conjecture}
Corollary~\ref{crochet en cyclique} follows directly from
Theorem~\ref{structure BV sur la cohomologie de Hochschild de M}
and from the following proposition.
In~\cite[Corollary 1.7 and Section 7]{MenichiL:BValgaccoHa},
we proved that if $A$ is a symmetric algebra then its negative cyclic
cohomology $HC^*_-(A)$ is a graded Lie algebra of lower degree 2.
In fact, we proved more generally
\begin{proposition}
If the Hochschild cohomology of a (differential graded) algebra $A$,
$HH^*(A;A^\vee)$,
equipped with $H_*(B^\vee)$,
has a Batalin-Vilkovisky algebra structure of degree $-d$
then its negative cyclic
cohomology $HC^*_-(A)$ is a graded Lie algebra of lower degree 2-d.
\end{proposition}
\begin{proof}
Apply \cite[Proposition 7.1]{MenichiL:BValgaccoHa} to the mixed complex
$\mathcal{C}^*(A;A^\vee)$ 
(desuspended $d$-times in order to take into account the degree $d$ shift).
By definition, $HC^*_-(A)$ is the differential torsion product

$
\hfill\text{Tor}^{H_*(S^1)}(\mathcal{C}^*(A;A^\vee),\mathbb{F}).
$
\end{proof}
Another interesting particular case
of~\cite[Proposition 7.1]{MenichiL:BValgaccoHa}
is the following proposition.
\begin{proposition}
If the Hochschild homology of an algebra $A$, $HH_*(A;A)$,
equipped with Connes boundary map $B$, has a
Batalin-Vilkovisky algebra structure then its cyclic
homology $HC_*(A)$ is a graded Lie algebra of lower degree 2.
\end{proposition}
\begin{proof}
Apply \cite[Proposition 7.1]{MenichiL:BValgaccoHa} to the mixed complex
$\mathcal{C}_*(A;A)$. By definition, $HC_*(A)$ is the differential
torsion product

$
\hfill\text{Tor}^{H_*(S^1)}(\mathcal{C}_*(A;A),\mathbb{F}).
$
\end{proof}
Remark that in fact, these graded Lie algebra structures extend
to $Lie_\infty$-algebra structures like
the Chas-Sullivan string bracket~\cite[Theorem 6.2 and Corollary 6.3]{Chas-Sullivan:stringtop}.

Chas-Sullivan string bracket is defined using Gysin long exact sequence.
The bracket given by Corollary~\ref{crochet en cyclique}
is defined similarly using Connes long exact sequence.
Jones~\cite{JonesJ:Cycheh} proved that Gysin and Connes long exact
sequences are isomorphic.
Therefore Conjecture~\ref{conjecture iso d'algebres} implies
Conjecture~\ref{conjecture iso d'algebres de Lie},
since as we explained in the introduction,
Conjecture~\ref{conjecture iso d'algebres} implies that
the Jones isomorphism
$$J:H_{p+d}(LM)\buildrel{\cong}\over\rightarrow
HH^{-p-d}(S^*(M),S^*(M)^\vee)$$
is an isomorphism of Batalin-Vilkovisky algebras.
\bibliography{Bibliographie}
\bibliographystyle{amsplain}
\end{document}